\documentclass[12pt,reqno]{amsart}

\usepackage{tabularx}
\usepackage{amssymb}        
\usepackage{latexsym}       
\usepackage{graphicx}       
\usepackage{multirow}       
\usepackage{xcolor}         
\usepackage[all]{xy}
\usepackage{caption,stmaryrd}
\usepackage{setspace}

\usepackage{enumitem}      
\setlist{
    topsep=6pt,           
    partopsep=0pt,        
    itemsep=4pt,          
    parsep=0pt,           
    leftmargin=*,         
    after=\addvspace{0.1cm}  
}
\setlist[enumerate]{font={\rm},itemsep=0.2\baselineskip}
\setlist[enumerate,1]{label={(\arabic*)}}
\setlist[enumerate,2]{label={(\roman*)}}

\usepackage{mathrsfs}
\usepackage{verbatim}
\usepackage{arydshln}
\usepackage{tikz}
\usetikzlibrary{positioning}

\usepackage[backref=page,hypertexnames=false]{hyperref}

\newtheorem{theorem}{Theorem}[section]

\newtheorem{lemma}[theorem]{Lemma}

\theoremstyle{definition}

\newtheorem{definition}[theorem]{Definition}

\newtheorem{problem}[theorem]{Problem}

\usepackage{etoolbox}
\newcommand{\declaremathfont}[2]{%
    \begingroup
    \newcommand{\declareone}[1]{%
        \expandafter\gdef\csname #1##1\endcsname{#2{##1}}%
    }%
    \forcsvlist{\declareone}{A,B,C,D,E,F,G,H,I,J,K,L,M,N,O,P,Q,R,S,T,U,V,W,X,Y,Z}%
    \endgroup
}
\declaremathfont{bb}{\mathbb}
\declaremathfont{cal}{\mathcal}
\declaremathfont{rm}{\mathrm}
\declaremathfont{bf}{\mathbf}
\declaremathfont{sf}{\mathsf}
\declaremathfont{tt}{\mathtt}
\declaremathfont{it}{\mathit}
\declaremathfont{frak}{\mathfrak}

\newcommand{\declaremathgroup}[1]{\expandafter\newcommand\csname #1\endcsname{{\mathrm{#1}}}}
\forcsvlist{\declaremathgroup}{%
    SL,GL,Sp,SU,SO,CSp,TL,
    PSL,PGL,PSp,PSU,
    ASL,AGL,
    Ree,Sz,Suz,Ru,Th,He,McL,HN,HS,Co,Ly,
}

\newcommand{\declaremathobject}[1]{\expandafter\newcommand\csname #1\endcsname{{\mathrm{#1}}}}
\forcsvlist{\declaremathobject}{%
    PG,GF,Gal,
    AF,HA,AS,PA,TW,SD, CD, HC,
    Af,Ha,As,Pa,Tw,Sd,Cd,
    KTS,STS,
}

\newcommand{\declaremathoperator}[1]{\expandafter\newcommand\csname #1\endcsname{{\mathsf{#1}}}}
\forcsvlist{\declaremathoperator}{%
    Ker,Tr,tr,rad,Rad,Soc,soc,Aut,Inn,Out,Hom,End,Syl,diag,Cay,Cos,Ann,rank,cl,core,Sym,
}

\def\GammaL{{\Gamma \mathrm{L}}}

\def\PGammaL{{\mathrm{P}\Gamma\mathrm{L}}}
\def\POmega{{\mathrm{P}\Omega}}

\def\CGammaSp{{\mathrm{C}\Gamma\mathrm{Sp}}}

\def\Magma{${\mathsf{MAGMA}}$}

\title{Finite groups with all involutions fused}

\author{Zheng Guo}
\address{School of Mathematics and Statistics, Zhengzhou University, Zhengzhou, Henan, China}
\email{12131227@mail.sustech.edu.cn {\text{\rm(Guo)}}}

\author{Yan Zhou Zhu}
\address{School of Mathematical Sciences, Xiamen University, Xiamen, Fujian, China}
\email{zhuyz@xmu.edu.cn {\text{\rm(Zhu)}}}

\date{}

\begin{document}

\begin{abstract}
We characterize finite groups whose involutions form a single fusion class.
In the solvable case, we show that the $2$-length is $1$, extending Thompson's theorem which states that a solvable group with all involutions conjugate has $2$-length $1$.
For a nonsolvable group $G$, we prove that if its involutions are fused then they are conjugate; moreover, $G/\bfO_{2'}(G)$ is almost quasisimple.
As an application, we obtain a restriction on $m$ for which an $m$-pyramidal Kirkman triple system can be realized under a nonsolvable group.
\end{abstract}


\maketitle

\section{Introduction}\label{introduction}

Let $G$ be a finite group.
Two elements $x,y\in G$ are said to be \textit{fused} if $x^{\sigma}=y$ for some $\sigma\in\Aut(G)$.
The orbits of $\Aut(G)$ on $G$ are called the \textit{fusion classes}, or \textit{automorphism orbits}, of $G$, and the same terminology applies to the subgroups of $G$.
Conjugacy and fusion are among the most fundamental topics in finite group theory, see~\cite{bors2023Automorphism,feit1989finite,hall1959constructions,zhang1992AT} and the references therein.
A natural problem is therefore to study finite groups in which a prescribed set of elements or of subgroups forms a single fusion class, which leads to the following problem.
\begin{problem}
    Characterize finite groups in which all elements (or subgroups) of order $p$, with $p$ prime, are conjugate or fused.
\end{problem}

The study of finite $p$-groups whose elements of order $p$ are fused has a long history.
It goes back to 1962, when Boen, Rothaus and Thompson~\cite{boen1962pautomorphic,boen1962Further} called a $p$-group \textit{$p$-automorphic} if its elements of order $p$ are fused.
For odd $p$, Shult proved that $p$-automorphic $p$-groups of odd order are exactly the homocyclic $p$-groups~\cite{shult1969odd}.
For $p=2$, however, the problem is far more subtle.
Higman introduced and classified the \textit{Suzuki $2$-groups}~\cite{higman1963Suzuki}, and Shaw showed that a nonabelian $2$-automorphic $2$-group with more than one involution and a solvable automorphism group is a Suzuki $2$-group~\cite{shaw1970Sylow}.
Gross conjectured that every nonabelian $2$-automorphic $2$-group with more than one involution is a Suzuki $2$-group~\cite{gross1976automorphic}, and after the partial results of Bryukhanova~\cite{bryukhanova1981Automorphism} and Wilkens~\cite{wilkens1996note}, this was proved in 2024 by Li and Zhu~\cite{li2024Gross}.
Thus the $p$-automorphic $p$-groups are precisely the homocyclic $p$-groups, the generalized quaternion $2$-groups and the Suzuki $2$-groups.

The main result of this paper is a characterization of finite groups whose involutions are fused.
It is worth noting that Thompson proved that if all involutions of a solvable group $G$ are conjugate, then $G$ has $2$-length $1$, that is, $G/\bfO_{2'}(G)$ has a normal Sylow $2$-subgroup~\cite[Lemma~5.40]{Thompson1968}.
We extend Thompson's theorem from conjugacy to fusion: a finite solvable group with fused involutions has $2$-length $1$.
For a nonsolvable group $G$, we prove that the involutions of $G$ are fused if and only if they are conjugate.
Moreover, if $G$ is nonsolvable and its involutions are fused, then $G/\bfO_{2'}(G)$ is an almost quasisimple group.
Yamaki classified finite simple groups with a unique conjugacy class of involutions~\cite{yamaki2008invconj}.
Recently, Li and Zeng extended Yamaki's classification to the almost simple case as one of the main ingredients in their solution to Brauer's Problem 33~\cite{arxiv-sylow-2-subgroup}. 
We now state our main result.
\begin{theorem}\label{thm:main_inv_fuse}
    Let $G$ be a finite group whose involutions are fused.
    Then either
    \begin{enumerate}
        \item $G$ is solvable with $2$-length $\ell_2(G)=1$, that is 
        \[ G=\bfO_{2'}(G).(B{:}C),\]
        where $B$ is a Sylow $2$-subgroup of $G$ and $C$ has odd order; or
        \item involutions of $G$ are conjugate and $G/\bfO_{2'}(G)$ is an almost quasisimple group with $(G/\bfO_{2'}(G))^{(\infty)}$ given in Table~$\ref{tab:pyra}$, see Theorem~$\ref{thm:as_inv_fuse}$ and Lemma~$\ref{lem:uni_inv}$ for details.
    \end{enumerate}
\end{theorem}

\begin{table}[htb]
    \centering
    \caption{Quasisimple groups whose involutions are conjugate}\label{tab:pyra}
    \footnotesize
    \renewcommand{\arraystretch}{1.15}
    \begin{tabular}{@{\hspace{.4em}}l@{\hspace{.9em}}l@{\hspace{.9em}}l@{\hspace{.9em}}l@{\hspace{.4em}}}
        \noalign{\hrule \vskip 1.5pt}
        $G$ & $k$: number of involutions  & condition &  $3\mid k$ \\
        \noalign{\vskip 1.5pt\hrule\vskip 2pt}
        $2.\rmA_7$ & $1$ & -- & $\times$ \\
        $\SL_2(q)$ & $1$ & $q\geqslant 5$ odd & $\times$ \\
        $\PSL_2(q)$ & $q^2-1$ & $q\geqslant4$ even & $\checkmark$ \\
         & $q(q+1)/2$& $q\equiv 1\pmod 4$ & $q\equiv 5,9\pmod{12}$ \\
         & $q(q-1)/2$& $q\equiv 3\pmod 4$, $q> 3$ & $q\equiv 3,7\pmod{12}$ \\
        $\PSL_3(q)$ & $q^2(q^2+q+1)$ & $q$ odd & $q\equiv 1,3\pmod 6$ \\
         & $(q^3-1)(q+1)$ & $q$ even  & $\checkmark$ \\
        $\PSL_4(q)$ & $q^4(q^2+q+1)(q^2+1)/2$ & $q\equiv 5\pmod 8$ & $q\equiv 13\pmod{24}$ \\
        $\PSU_3(q)$ & $q^2(q^2-q+1)$ & $q$ odd  & $q\equiv 3,5\pmod 6$ \\
         & $(q^3+1)(q-1)$ & $q\geqslant4$ even & $\checkmark$ \\
        $\PSU_4(q)$ & $q^4(q^2-q+1)(q^2+1)/2$ & $q\equiv 3\pmod 8$ & $q\equiv 3,11\pmod{24}$ \\
        ${}^{2}\rmG_2(q)$ & $q^2(q^2-q+1)$ & $q=3^{2m+1}>3$ & $\checkmark$ \\
        ${}^{3}\rmD_4(q)$ & $q^8(q^8+q^4+1)$ & $q$ odd & $\checkmark$ \\
        $\rmG_2(q)$ & $q^4(q^4+q^2+1)$ & $q$ odd & $\checkmark$ \\
        ${}^{2}\rmB_2(q)$ & $(q^2+1)(q-1)$ & $q=2^{2m+1}>2$ & $\times$ \\
        \noalign{\vskip 1.5pt\hrule\vskip 2pt}
        $\rmJ_1$ & $1\,463$ & -- & $\times$ \\
        $\rmA_7$ & $105$ & -- & $\checkmark$ \\
        $\rmM_{11}$ & $165$ & -- & $\checkmark$ \\
        $\rmM_{22}$ & $1\,155$  & -- & $\checkmark$ \\
        $\rmM_{23}$ & $3\,795$  & -- & $\checkmark$ \\
        $\rmJ_3$ & $26\,163$  & -- & $\checkmark$ \\
        $\McL$ & $22\,275$  & -- & $\checkmark$ \\
        $\mathrm{O'N}$ & $2\,857\,239$  & -- & $\checkmark$ \\
        $\Th$ & $976\,841\,775$  & -- & $\checkmark$ \\
        $\Ly$ & $1\,296\,826\,875$  & -- & $\checkmark$ \\
        \noalign{\vskip 2pt\hrule}
    \end{tabular}
\end{table}

Another motivation for our study is from a problem regarding Steiner triple systems.
Recently, the groups that have a unique conjugacy class of involutions play an important role in the study of \textit{Steiner triple systems} and \textit{Kirkman triple systems}.
A Steiner triple system is a pair $(V,\calB)$, where $V$ is a set of $v$ points and $\calB$ is a set of $3$-subsets of $V$ such that any two distinct points lie in exactly one element of $\calB$, denoted by $\STS(v)$.
An $\STS(v)$ is called \textit{$m$-pyramidal} if its automorphism group has a subgroup $G$ which fixes $m$ points and acts regularly on the other $v-m$ points.
Buratti, Rinaldi and Traetta~\cite{buratti2017threepyramidal} observed that a nontrivial group $G$ realizing an $m$-pyramidal $\STS(v)$ has exactly $m$ involutions, and they proved that a $3$-pyramidal $\STS(v)$ exists if and only if $v\equiv7,9,15\pmod{24}$ or $v\equiv3,19\pmod{48}$.
A Kirkman triple system is a Steiner triple system together with a resolution of $\calB$ into parallel classes, denoted by $\KTS(v)$.
This leads to the following natural problem.

\begin{problem}
    For which pair $(v,m)$ does an $m$-pyramidal $\KTS(v)$ exist.
\end{problem}

Gao and Garonzi proved that if $G\leqslant\Aut(\KTS(v))$ fixes $m$ points and acts regularly on the other $v-m$ points, then $G$ has precisely $m$ involutions and these $m$ involutions are pairwise conjugate, see~\cite[Proposition 1.1]{gao2025}.
Accordingly, a finite group is called \textit{$m$-pyramidal} if it has precisely $m$ involutions and these $m$ involutions are pairwise conjugate, and \textit{pyramidal} if it is $m$-pyramidal for some $m$.
For $m=3$ necessary conditions and infinite families whose orders fill a congruence class are given in~\cite{Simona2021}; and no example seems to be known for $m>3$, as the discussion in~\cite{gao2025}.
For more on the existence problem for $m$-pyramidal $\STS$ and $\KTS$, we refer the readers to
\cite{buratti2017threepyramidal,chang2025pyramidalAbelian,Simona2021} and the references therein.
Gao and Garonzi then proved that, for a prime power $m=p^k$ for $p\neq 7$, each $m$-pyramidal group is solvable if and only if both $p,k$ are odd or $m=9$~\cite{gao2025}.
Applying Theorem~\ref{thm:main_inv_fuse} to the nonsolvable case, we obtain the following restriction for $m$-pyramidal $\KTS(v)$ realized under a nonsolvable group.

\begin{theorem}\label{thm:KTS}
    Assume that an $m$-pyramidal $\KTS(v)$ is realized under a nonsolvable group $G$.
    Let $O=\bfO_{2'}(G)$.
    Then the following statements hold.
    \begin{enumerate}
        \item $G/O$ is an almost quasisimple $k$-pyramidal group with $(G/O)^{(\infty)}$ isomorphic to one of the groups in Table~$\ref{tab:pyra}$.
        \item $m/k=|O|/|\bfC_{O}(x)|$ is an odd square number for each involution $x\in G$.
        \item one of the following holds.
        \begin{enumerate}
            \item $m/k$ is divisible by $3$.
            \item $m/k$ is not divisible by $3$, and $3\mid k$.
            \item $(G/O)^{(\infty)}\cong {}^2\rmB_2(q)$ with $q=2^{2n+1}$ for $n\geqslant 1$ and $O$ is a $\{2,3\}'$-group. 
        \end{enumerate}
    \end{enumerate}
    In particular, for a prime power $m$, each $m$-pyramidal group is solvable if and only if $m=p^{2k+1}$ for some odd prime $p$ or $m=9$.
\end{theorem}

\subsection*{Terminologies and notation}
Throughout, $G$ denotes a finite group.
For a prime number $p$, denote by $\bfO_{p}(G)$ (or $\bfO_{p'}(G)$, resp.) the maximal normal $p$-subgroup (or $p'$-subgroup, resp.) of $G$.
We use $G^{(\infty)}$ for the last term in the derived series of $G$.
We say $G$ is \textit{almost simple} if $\soc(G)$ is a nonabelian simple group, and say $G$ is \textit{quasisimple} if $G$ is a perfect group (i.e. $G'=G$) such that $G/\bfZ(G)$ is a nonabelian simple group.
If $G^{(\infty)}$ is quasisimple and $G/\bfZ(G^{(\infty)})$ is almost simple, then we call $G$ an \textit{almost quasisimple} group.

\section{Preliminaries}\label{sec:pre}

In this section, we collect some elementary facts and some known results.
We start with some basic facts in the following lemma.

\begin{lemma}\label{lem:pre}
    Let $G$ be a finite group and let $N$ be a characteristic subgroup of $G$.
    Then the following statements hold.
    \begin{enumerate}
        \item If $x,y\in N$ are fused in $G$, then they are fused in $N$.
        \item If $x,y\in G\setminus N$ are fused in $G$, then $xN$ and $yN$ are fused in $G/N$.
        \item Assume that $N=R_1\times \cdots\times R_n$ for quasisimple subgroups $R_1,...,R_n$ with $n\geqslant 2$.
        Then, for any $1\leqslant i<j\leqslant n$ and any $1\neq r_i\in R_i$, the elements $r_i$ and $r_ir_j$ are not fused.
    \end{enumerate}
\end{lemma}
\begin{proof}
    Parts~(1) and (2) are clear as every automorphism of $G$ induces an automorphism of $N$ and $G/N$.

    Note that there is no quasisimple normal subgroup of $N$ containing $r_ir_j$.
    Since $R_i$ is a quasisimple normal subgroup of $N$, we have that $r_i$ and $r_ir_j$ are not fused in $N$.
    By part~(1) $r_i$ and $r_ir_j$ are not fused in $G$, as in part~(3).
\end{proof}

A finite $p$-group with elements of order $p$ fused is called $p$-automorphic.
For odd prime $p$, Shult proved that a $p$-automorphic $p$-group of odd order is abelian~\cite[Corollary 2]{shult1969odd}, which forces it to be homocyclic.
The situation for $p=2$ is much richer, and the study goes back to G.~Higman~\cite{higman1963Suzuki}, who introduced and classified the \textit{Suzuki $2$-groups}.
\begin{definition}
    Suzuki $2$-groups are nonabelian $2$-groups $G$ with more than one involution which admit a cyclic subgroup $C\leqslant\Aut(G)$ acting transitively on the involutions of $G$.
\end{definition}

Suzuki $2$-groups are special $2$-groups of order $2^{n+n}$ ($n\geqslant 3$) or $2^{n+2n}$ ($n\geqslant 3$), where $p^{a+b}$ stands for a special $p$-group of order $p^{a+b}$ with Frattini subgroup of order $p^a$.
Building on the work of Shaw~\cite{shaw1970Sylow}, Gross~\cite{gross1976automorphic}, Bryukhanova~\cite{bryukhanova1981Automorphism} and Wilkens~\cite{wilkens1996note}, Li and Zhu~\cite{li2024Gross} completed the classification of $p$-automorphic $p$-groups for all $p$, as in the following theorem.
\begin{theorem}\label{thm:p-group}
    Let $P$ be a finite $p$-group such that $\Aut(P)$ is transitive on the set of all elements of order $p$. 
    Then one of the following statements holds.
    \begin{enumerate}
        \item $P$ is a homocyclic $p$-group;
        \item $P$ is a generalized quaternion $2$-group;
        \item $P$ is a Suzuki $2$-group and $\Aut(P)$ is solvable.
    \end{enumerate}
\end{theorem}

Yamaki determined the finite simple groups with a unique conjugacy class of involutions and proved that if a finite simple group has more than one conjugacy class of involutions, then it has more than one fusion class of involutions, see \cite[Lemmas 1 \& 2]{yamaki2008invconj}.
Recently, Li and Zeng extended Yamaki's result by determining the finite almost simple groups with a unique conjugacy class of involutions~\cite{arxiv-sylow-2-subgroup}.
Consequently, the groups in Li-Zeng's list are precisely the finite almost simple groups with a unique fusion class of involutions, as stated in the following theorem.
For a group of Lie type $G$, we denote by $\delta$ and $\phi$ a diagonal automorphism and a field automorphism of $G$, respectively.
\begin{theorem}\label{thm:as_inv_fuse}
    Let $G$ be a finite almost simple group whose involutions are fused.
    Then $G$ has a unique conjugacy class of involutions, and $G$ equals one of the following groups.
    \begin{enumerate}
        \item $\rmA_7$, $\rmM_{11}$, $\rmM_{22}$, $\rmM_{23}$, $\rmJ_1$, $\rmJ_3$, $\McL$, $\mathrm{O'N}$,\  $\Ly$,\ $\Th$,
        \item $\PSL_2(p^f)\trianglelefteqslant G\leqslant\PSL_2(p^f).(\langle\delta\phi^{\frac{f}{2}}\rangle\times \langle \phi\rangle_{2'})$ for odd prime $p$ and even $f$, and 
        \item $S\trianglelefteqslant G\leqslant S{:}\Out(S)_{2'}$, where $S$ is isomorphic to one of the following groups:
        \[
        \begin{array}{l}
        \PSL_2(q),\ \PSL_3(q),\ \PSL_4(q)\ (q\equiv 5\ \mathrm{mod}\ 8),\ \PSU_3(q),\ \PSU_4(q)\ (q\equiv 3\ \mathrm{mod}\ 8),\\[2pt]
        {}^2\rmG_2(3^{2m+1}),\ {}^2\rmB_2(2^{2m+1}),\ {}^3\rmD_4(q)\ (q\ \mathrm{odd})\mbox{ and } \rmG_2(q)\ (q\ \mathrm{odd}).
        \end{array}
        \]
    \end{enumerate}
\end{theorem}

It is well-known that a $p$-group with a unique subgroup of order $p$ is cyclic or generalized quaternion, see~\cite[page 199]{gorenstein1980Finite} for instance. 
Hence finite groups with a unique involution have cyclic or generalized quaternion Sylow $2$-subgroups.
Finite groups with cyclic Sylow $2$-subgroups are solvable, and finite groups with dihedral Sylow $2$-subgroups are well-characterized by Gorenstein and Walter~\cite{gorenstein1965SylDih}.
Combining these results, we obtain the following characterization of finite groups with a unique involution.
Before that, we define 
\[\TL_2(q)=\langle M\in\SL_2(q^2)\mid M^{\varphi}=M\text{ or }M^{\varphi}=-M\rangle,\]
where $\varphi$ maps each entry of $M$ to its $q$-th power.
For odd $q$, one can check that $\SL_2(q)\lhd \TL_2(q)$ has index $2$, and $\TL_2(q)\cong 2.\PGL_2(q)\cong \SL_2(q).\langle \delta\rangle$, where $\delta$ is a diagonal automorphism of $\SL_2(q)$.
\begin{lemma}\label{lem:uni_inv}
    Assume that $G$ is a finite nonsolvable group with a unique involution.
    If $\bfO_{2'}(G)=1$, then $G$ is an almost quasisimple group.
    Moreover, one of the following holds.
    \begin{enumerate}
        \item $G\cong 2.\rmA_7$.
        \item $\SL_2(q)\trianglelefteqslant G\leqslant \SL_2(q){:}\langle\phi\rangle_{2'}$ for odd $q\geqslant 5$.
        \item $\TL_2(q)\trianglelefteqslant G\leqslant \TL_2(q){:}\langle\phi\rangle_{2'}$ for odd $q\geqslant 5$.
    \end{enumerate}
\end{lemma}
\begin{proof}
    Let $P$ be a Sylow $2$-subgroup of $G$, and let $x$ be the unique involution in $G$.
    Then $x$ is the unique involution of $P$.
    Hence $P$ is either cyclic or generalized quaternion.
    If $P$ is cyclic, then $G$ is solvable which is a contradiction.
    Thus $P\cong \rmQ_{2^k}$ with $k\geqslant 3$.
    Note that $G/\langle x\rangle$ has a Sylow $2$-subgroup $P/\langle x\rangle\cong \rmD_{2^{k-1}}$.
    Since the full preimage of $\bfO_{2'}(G/\langle x\rangle)$ in $G$ is $\langle x\rangle \times \bfO_{2'}(G)=\langle x\rangle$, we have that $\bfO_{2'}(G/\langle x\rangle)=1$.
    Hence, by~\cite[Theorem 1]{gorenstein1965SylDih}, either 
    \[G/\langle x\rangle\lesssim \PGammaL_2(q)\mbox{ with $(G/\langle x\rangle)^{(\infty)}\cong\PSL_2(q)$ for $q\geqslant 4$, or }G/\langle x\rangle\cong\rmA_7.\]
    Since $\langle x\rangle\leqslant \bfZ(G)$ and $P$ is non-split over $\langle x\rangle$, we have that $G^{(\infty)}$ is a quasisimple group.
    This implies that 
    \[G\cong 2.\rmA_7\mbox{ or }\SL_2(q)\trianglelefteqslant G\mbox{ for $q\geqslant 4$}.\]
    Assume the latter case holds, that is, $\SL_2(q)\trianglelefteqslant G$. 
    Note that $G/\langle x\rangle\leq\PGammaL_2(q)$ has a dihedral Sylow $2$-subgroup.
    By~\cite[page 510]{suzuki1986Group}, we have that either 
    \[G/\langle x\rangle\cong \PSL_2(q){:}\langle\phi\rangle_{2'}\mbox{ or }G/\langle x\rangle\cong\PGL_2(q){:}\langle\phi\rangle_{2'},\]
    where $\phi$ is a field automorphism of $\PSL_2(q)$.
    Note that each of $\SL_2(q)$ and $TL_2(q)$ has a unique involution.
    By~\cite[Theorem 4.1]{babai2000Automorphisms}, we have that either
    \[G\cong \SL_2(q){:}\langle\phi\rangle_{2'}\mbox{ or }G\cong\TL_2(q){:}\langle\phi\rangle_{2'}.\]
    Hence $G$ satisfies one of the three cases.
    This completes the proof of this lemma.
\end{proof}

We end this part with the following lemma on almost quasisimple groups.
Note that the \textit{Schur multiplier} of a simple group $S$ is the center of the largest quasisimple group $G$ with $G/\bfZ(G)\cong S$.

\begin{lemma}\label{lem:aqs_inv}
    Let $G$ be an almost quasisimple group with $\bfO_2(G)\neq 1$ and $\bfO_{2'}(G)=1$. 
    If all involutions in $G$ are fused, then $G$ contains a unique involution and $G$ is one of the groups given in Lemma~\ref{lem:uni_inv}.
\end{lemma}
\begin{proof}
    Note that every involution of $G$ lies in $\bfO_2(G)\leqslant\bfZ(G^{(\infty)})$.
    If $\bfZ(G^{(\infty)})$ contains exactly one involution, then $G$ has a unique involution.
    Then the lemma holds in this case by Lemma~\ref{lem:uni_inv}.
    Suppose that $\bfZ(G^{(\infty)})$ has more than one involution.
    By the classification of finite simple groups together with their Schur multiplier, we have that $\bfZ(G^{(\infty)})$ has exactly $3$ involutions.
    As $\Aut(G)$ is transitive on the three involutions of $G$, there exists an element $\sigma\in\Aut(G)$ of order $3^k$ for some $k\geqslant 1$ such that $\langle\sigma\rangle$ is transitive on involutions of $G$.
    Define $X=G{:}\langle \sigma\rangle$.
    Then $X$ has exactly $3$ involutions and the involutions of $X$ are conjugate.
    By~\cite{gao2023threepyramidal}, $X$ is solvable, a contradiction.
    Therefore, the proof is complete.
\end{proof}

\section{Finite groups with unique fusion classes of involutions}


\subsection{A reduction theorem}\

In this subsection, we aim to prove the following reduction lemma which is a primary structure theorem for finite groups with one fusion class of involutions.
\begin{theorem}\label{thm:reduction}
    Let $G$ be a finite group with a unique fusion class of involutions.
    Assume that $\bfO_{2'}(G)=1$.
    Then one of the following statements holds.
    \begin{enumerate}
        \item $G$ is a solvable group;
        \item $G$ is an almost quasisimple group;
        \item $G$ has a unique conjugacy class of involutions such that $\bfO_2(G)\cong \rmC_{2^m}^n$ is homocyclic and is self-centralized in $G$.
    \end{enumerate}
\end{theorem}

In order to prove the above theorem, we consider the \textit{generalized Fitting subgroup} $\rmF^*(G)$ of a finite group $G$.
The maximal nilpotent normal subgroup of $G$ is called the \textit{Fitting subgroup} of $G$, denoted by $\rmF(G)$; and the \textit{layer} of $G$ is generated by all quasisimple subnormal subgroups of $G$, denoted by $\rmE(G)$.
Then 
\[\rmF^*(G):=\rmF(G)\rmE(G) \text{ and }\rmE(G)=(R_1\times \cdots\times R_n)/N,\] 
where $R_1,...,R_n$ are quasisimple groups and $N\leqslant \bfZ(R_1)\times\cdots\times \bfZ(R_n)$.
Note that $[\rmF(G),\rmE(G)]=1$ and $\rmF^*(G)$ is self-centralized.
We refer readers to~\cite[Chapter 11]{aschbacher2000Finite} for more properties of the generalized Fitting subgroup.

\begin{lemma}\label{lem:reduction-E}
    Let $G$ be a finite group with a unique fusion class of involutions.
    Assume that $\bfO_{2'}(G)=1$ and $\rmE(G)\neq 1$.
    Then $G$ is an almost quasisimple group.
\end{lemma}
\begin{proof}
    Note that $\rmE(G)\neq 1$ is of even order.
    Hence every involution of $G$ lies in $\rmE(G)$.
    As $\bfO_{2'}(G)=1$, either $\rmE(G)=T_1\times\cdots\times T_n$ for nonabelian simple groups $T_1,...,T_n$ or $\bfZ(\rmE(G))=\bfO_2(G)\neq 1$.

    Assume that $\rmE(G)=T_1\times\cdots\times T_n$ for nonabelian simple groups $T_1,...,T_n$.
    Let $t_i$ be an involution in $T_i$ for $i=1,...,n$.
    Lemma~\ref{lem:pre}\,(3) implies that $t_1$ and $t_1t_2$ are not fused if $n\geqslant 2$.
    Thus $n=1$, and then $\rmE(G)$ is a nonabelian simple group.
    Hence $\bfO_2(G)=1$ and $\rmF^*(G)=\rmF(G)\times \rmE(G)$.
    Since $\rmF(G)$ is a direct product of $\bfO_p(G)$ for all prime divisors $p$ of $|G|$, we have that $\rmF(G)=1$.
    Thus $\rmF^*(G)=\rmE(G)$, and then $G\lesssim \Aut(\rmE(G))$ is an almost simple group.

    Assume that $\bfZ(\rmE(G))=\bfO_2(G)\neq 1$ contains all involutions of $G$.
    Then involutions of $G$ are all contained in $\rmE(G)$ and they are fused in $\rmE(G)$ by Lemma~\ref{lem:pre}\,(1).
    Set 
    \[\rmE(G)=(R_1\times \cdots\times R_n)/N,\] 
    where $R_1,...,R_n$ are all quasisimple groups and $N\leqslant \bfZ(R_1)\times\cdots\times \bfZ(R_n)$.
    Denote by $\overline{\cdot}$ the quotient map from $R_1\times\cdots\times R_n$ to $\rmE(G)$, and let $Z=\overline{\bfZ(R_1)\times\cdots\times \bfZ(R_n)}=\bfZ(\rmE(G))=\bfO_2(G)$.
    We deal with this case step by step.

    \vspace{.2cm} 
    \noindent\textbf{Step 1.} $\bfZ(R_i)$ is cyclic for each $i=1,...,n$.

    Note that $\bfZ(R_i)=\bfO_2(R_i)$.
    We may suppose that $\bfO_2(R_1)$ is not cyclic.
    Then the Schur multiplier of $S=R_1/\bfZ(R_1)$ has more than one involution.
    By the classification of finite simple groups, we summarize that if the Schur multiplier of the simple group $S$ has more than one involution, then $S$ is one of the following groups whose $2$-multiplier has abelian rank $2$:
    \[\PSL_3(4),\ \POmega_{2n}^+(q)\mbox{ (even $n\geqslant 4$ and odd $q$)},\ \POmega_8^+(2),\ \PSU_6(2),\ {}^2\rmE_6(2^2)\mbox{ and }{}^2\rmB_2(8).\]
    Suppose that $R_1$ contains more than $3$ involutions. 
    Then there exists an involution $x\in R_1\setminus \bfZ(R_1)$.
    Hence $x N$ is an involution in $\rmE(G)\setminus \bfZ(\rmE(G))$, a contradiction.
    Thus we only need to show that $R_1$ has more than $3$ involutions case by case.
    
    \textbf{Case 1.} Assume that $S\cong \POmega_{8}^+(q)$ for odd $q$.
    
    Then $R_1\cong \rmC_2^2.S$ and there exists an automorphism $\sigma\in\Aut(S)$ of order $3$ that permutes the three involution in $\bfZ(R_1)$ transitively, see~\cite[Section 4.7.2]{wilson2009finite}.
    Hence $S$ has more than $3$ involutions, for otherwise $R{:}\langle \sigma\rangle$ is a nonsolvable group which contains exactly three involutions and all three involutions are conjugate.
    However, a finite group containing exactly $3$ involutions and a unique conjugacy class of involutions is solvable by~\cite{gao2023threepyramidal}, a contradiction.

    \textbf{Case 2.} Assume that $S\cong \PSL_3(4)$ or ${}^2\rmB_2(8)$.

    Assume that $S\cong \PSL_3(4)$.
    Then $\bfZ(R_1)$ a quotient of $\rmC_4^2$.
    Hence $R_1$ is a quotient of the quasisimple group $R:=\rmC_4^2.\PSL_3(4)$, that is, $R_1=R/K$ with $K\leqslant\bfZ(R)$.
    With the aid of \Magma, we have that $R$ contains an involution $x$ in $R\setminus\bfZ(R)$.
    Then $xK$ is an involution of $R_1\setminus\bfZ(R_1)$.
    Recall that $\bfZ(R_1)$ is not cyclic.
    Thus $R_1$ has more than $3$ involutions.
    
    Assume that $S\cong {}^2\rmB_2(8)$.
    Then $R_1\cong \rmC_2^2.S$.
    The permutation representation of $R_1{:}\rmC_3$ of degree $2080$ is available in~\cite{webatlas}.
    Then, with the aid of \Magma~\cite{magma}, we have that $\bfZ(R_1{:}\rmC_3)=1$.
    This implies that the three involutions in $\bfZ(R_1)$ are conjugate in $R_1{:}\rmC_3$.
    The same argument as in Case 1 yields a contradiction.

    \textbf{Case 3.} Assume that $\POmega_{2n}^+(q)$ (even $n\geqslant 6$ and odd $q$), $\POmega_{8}^+(2)$, $\PSU_6(2)$ or ${}^2\rmE_6(2^2)$.

    We claim that $S$ contains a subgroup $\overline{P}$ isomorphic to $\rmC_2^5$:
    the maximal torus of $\POmega_{2n}^+(q)$ contains $\rmC_2^n$ for $n\geqslant 6$ and odd $q$;
    ${}^2\rmE_6(2^2)$ has a subgroup isomorphic to $\rmC_2^8$ by~\cite{ATLAS};
    for the other cases, we can verify the claim with the aid of \Magma.
    Let $P$ be the full preimage of $\overline{P}$ in $R_1$.
    Then $P\cong \rmC_2^2.\rmC_2^5$ has exponent $4$ and has exactly $3$ involutions.
    By verifying all groups of order $2^7$ and exponent $4$ in \Magma, we have that every such group has more than $3$ involutions.
    Thus $R_1$ has more than $3$ involutions.

    \vspace{.2cm} 
    \noindent\textbf{Step 2.} $R_i\cong\overline{R_i}\cong \SL_2(q)$ for some odd $q$ or $2.\rmA_7$ for each $i=1,...,n$.

    As $\bfZ(R_i)$ is cyclic, $\bfZ(\overline{R_i})$ is also cyclic.
    Note that involutions of $\overline{R_i}$ are contained in $Z\cap\overline{R_i}=\bfZ(\overline{R_i})$.
    Thus $\bfZ(\overline{R_i})\neq 1$ and $\overline{R_i}$ has a unique involution.
    By Lemma~\ref{lem:uni_inv}, we have that $\overline{R_i}\cong \SL_2(q)$ for some odd $q$ or $2.\rmA_7$, so does $R_i$.

    \vspace{.2cm} 
    \noindent\textbf{Step 3.} $N=1$ and $\rmE(G)=R_1\times\cdots \times R_n$.

    Suppose that $z\in N$ is a nontrivial involution.
    Then $z=(r_1,...,r_n)\in R_1\times\cdots\times R_n$ such that $r_i\in\bfZ(R_i)$.
    We may assume that $r_1,...,r_k\neq 1$ and $r_{k+1},...,r_n=1$.
    By Step 2, $R_i\cong \SL_2(q)$ for some odd $q$ or $2.\rmA_7$.
    Then $r_i$ is an involution in $R_i$ and there exists $t_i\in R_i\setminus\bfZ(R_i)$ of order $4$ such that $t_i^2=r_i$ for each $i=1,...,k$.
    Let $x=(t_1,...,t_k,1,...,1)\in R_1\times \cdots\times R_n$.
    Then $x^2=z\in N$, and hence $xN$ is an involution in $\rmE(G)$.
    Note that $x\notin \bfZ(R_1)\times\cdots\times \bfZ(R_n)$.
    Thus $xN\notin Z$, and hence $G$ contains an involution outside $Z$, a contradiction.

    \vspace{.2cm} 
    \noindent\textbf{Step 4.} $G$ is an almost quasisimple group.

    Note that $\rmE(G)=R_1\times\cdots\times R_n$ is a direct product of quasisimple groups.
    By Lemma~\ref{lem:pre}\,(3), involutions of $\rmE(G)$ are not fused in $\rmE(G)$ if $n\geqslant 2$.
    So $n=1$ and $\rmE(G)$ is a quasisimple group.
    By Step 2, $\rmE(G)\cong \SL_2(q)$ for some odd $q$ or $2.\rmA_7$.
    Then $\bfZ(\rmE(G))\cong \rmC_2$.
    Since $\bfZ(\rmE(G))$ is a characteristic subgroup of $G$, $G$ has a unique involution.
    By Lemma~\ref{lem:uni_inv}, we have that $G$ is an almost quasisimple group.
\end{proof}

Now we consider the case where the generalized Fitting subgroup equals the Fitting subgroup.
\begin{lemma}\label{lem:reduction-F}
    Let $G$ be a finite group with a unique fusion class of involutions.
    Assume that $\bfO_{2'}(G)=\rmE(G)=1$.
    Then $\rmF^*(G)=\rmF(G)=\bfO_2(G)$, and either 
    \begin{enumerate}
        \item $G$ is solvable; or
        \item $G$ has a unique conjugacy class of involutions, $\bfO_2(G)\cong \rmC_{2^m}^n$ and $G/\bfO_2(G)\lesssim \GL_n(2)$ is faithful and transitive on non-identities of both 
        \[\Omega_1(\bfO_2(G))=\{x\in\bfO_2(G):x^2=1\}\mbox{ and }\bfO_2(G)/\Phi(\bfO_2(G)).\]
    \end{enumerate} 
\end{lemma}
\begin{proof}
    If $G$ is solvable, then case~(1) holds and there is nothing to prove.
    Hence we assume throughout the remainder of the proof that $G$ is nonsolvable.
    
    Note that $\rmF^*(G)=\rmF(G)$ is a product of $\bfO_p(G)$ for prime divisors $p$ of $|G|$.
    Thus $\rmF^*(G)=\rmF(G)=\bfO_2(G)$, and hence it is a $2$-automorphic $2$-group by Lemma~\ref{lem:pre}\,(1).
    By Theorem~\ref{thm:p-group}, $\rmF^*(G)$ equals $\rmC_{2^m}^n$ for $m,n\geqslant 1$, or $\rmQ_{2^k}$ for $k\geqslant 3$, or a Suzuki $2$-group with $\Aut(\rmF^*(G))$ solvable.
    Then $G/\rmF^*(G)\lesssim\Out(\rmF^*(G))$ is nonsolvable, and hence $\rmF^*(G)=\rmF(G)=\bfO_2(G)\cong \rmC_{2^m}^n$ is homocyclic and self-centralized in $G$.
    
    Let $N=\bfO_2(G)$.
    Then there is a natural group homomorphism $\psi:\Aut(G)\rightarrow \Aut(N/\Phi(N))$.
    Since $\Omega_1(N)=\{x\in N:x^2=1\}$ consists of the identity and all involutions of $G$, $\Aut(G)$ acts transitively on non-identities of $\Omega_1(N)=\{x\in N:x^2=1\}\cong \rmC_2^n$.
    The actions of $\Aut(G)$ on $\Omega_1(N)$ and $N/\Phi(N)$ are equivalent, and hence $G$ acts transitively on the non-identities of $N/\Phi(N)$.
    This implies that the image $\psi(\Aut(G))$ is a transitive linear group on $N/\Phi(N)\cong \bbF_2^n$.
    Note that $\Inn(G)\cong G$ is nonsolvable and the kernel of $\psi|_{\Inn(G)}$ is $N$.
    Thus $\psi(\Inn(G))$ is a nonsolvable normal subgroup of $\psi(\Aut(G))$.
    Transitive linear groups are stabilizers of affine $2$-transitive groups, and are classified by Huppert and Hering~\cite{hering1985Transitive,huppert1957Zweifach}.
    From the classification, we have that every nonsolvable normal subgroup of a nonsolvable transitive linear group on $\bbF_2^n$ is also a transitive linear group.
    Thus $\psi(\Inn(G))$ is transitive on the non-identities of $N/\Phi(N)$.
    This implies that $G$ is transitive on the non-identities of both $N/\Phi(N)$ and $\Omega_1(N)$.
    Hence involutions of $G$ are conjugate.

    It remains to show that $G/N$ is faithful on $\Omega_1(N)$.
    Note that the kernel of $\Aut(N)$ acting on $N/\Phi(N)$ is a $2$-group, see~\cite[Theorem 3.17]{huppert2025Finite} for instance.
    Since $\bfO_2(G/N)=1$, we have that $G/N\lesssim\Aut(N/\Phi(N))\cong\GL_n(2)$ is faithful on $N/\Phi(N)$.
    Then $G/N$ is faithful on non-identities of both $\Omega_1(N)$ and $N/\Phi(N)$.
\end{proof}

Now we summarize the proof of Theorem~\ref{thm:reduction}.
\begin{proof}[Proof of Theorem~\ref{thm:reduction}.]
    Since $\bfO_{2'}(G)=1$, the Fitting subgroup $\rmF(G)$ is a $2$-group, that is, $\rmF(G)=\bfO_2(G)$, and hence $\rmF^*(G)=\bfO_2(G)\rmE(G)$ such that $[\bfO_2(G),\rmE(G)]=1$.
    If $\rmE(G)\neq 1$, then $G$ is almost quasisimple with $G^{(\infty)}\cong\SL_2(q)$ for some odd $q$ or $2.\rmA_7$ by Lemma~\ref{lem:reduction-E}, which is case~(2).
    If $\rmE(G)=1$, then $\rmF^*(G)=\rmF(G)=\bfO_2(G)$.
    Lemma~\ref{lem:reduction-F} shows that $G$ is either solvable, which is case~(1), or has a unique conjugacy class of involutions with $\bfO_2(G)\cong\rmC_{2^m}^n$ homocyclic and self-centralized in $G$, which is case~(3).
    Thus one of the three statements holds.
\end{proof}

\subsection{Solvable groups whose involutions are fused}~\label{sov-2-fused}\

Thompson proved that if a solvable group has a non-cyclic Sylow
$2$-subgroup and all its involutions are conjugate, then it has $2$-length
$1$~\cite{Thompson1968}.
We will prove that a finite solvable group with a unique fusion class of
involutions also has $2$-length $1$, as stated in Theorem~\ref{thm:sol_inv_fuse}.
This generalizes Thompson's theorem.

\begin{theorem}\label{thm:sol_inv_fuse}
    Assume that $G$ is a solvable group whose involutions are fused.
    Then 
    \[G=A.(B{:}C)=(A{:}B).C,\]
    where both $A$ and $C$ have odd order and $B$ is a group stated in Theorem~\ref{thm:p-group}.
\end{theorem}

Before proving this theorem, we present a simple lemma.
\begin{lemma}\label{lem:splitnormal}
    Let $G$ be a finite group with a normal subgroup $A{:}B$ such that $\gcd(|A|,|B|)=1$ and $\bfC_A(B)=1$.
    Then $G=A{:}\bfN_G(B)$.
\end{lemma}
\begin{proof}
    For any $g\in G$, $B^g< A{:}B$ and $|B^g|$ is coprime to $|A|$. 
    Hence we have that $B^g$ is a complement of $A$ in $A{:}B$.
    Then there exists $x\in A$ such that $B^g=B^x$ by the Schur-Zassenhaus Theorem as $\gcd(|A|,|B|)=1$.
    It follows that $gx^{-1}$ normalizes $B$.
    Since $\bfN_A(B)=\bfC_A(B)=1$, we have that $B^y\neq B$ for any nontrivial element $y\in A$.
    Note that $gx^{-1}\in gA$.
    Then, for each coset $C$ in $G/A$, there exists an element in $C$ normalizing $B$.
    This implies that $|\bfN_G(B)|\geqslant |G|/|A|$, and hence $|G|\leqslant |A||\bfN_G(B)|$.
    Since $\bfN_G(B)\cap A=\bfN_A(B)=\bfC_A(B)=1$, we have that 
    \[|A\bfN_G(B)|=|A||\bfN_G(B)|\geqslant |G|.\]
    Thus $G=A\bfN_G(B)$.
    As $A$ is a normal Hall subgroup of $A{:}B$, $A$ is also normal in $G$.
    Therefore, $G=A{:}\bfN_G(B)$.
\end{proof}

Now we are ready to prove Theorem~\ref{thm:sol_inv_fuse}.
\begin{proof}[Proof of Theorem~\ref{thm:sol_inv_fuse}]
    Let $A=\bfO_{2'}(G)$, and let $X=G/A$.
    By Lemma~\ref{lem:pre}\,(2), involutions of $X$ are fused.
    Set $\rmF$ the Fitting subgroup of $X$.
    Then $\rmF=\bfO_2(X)$ as $\bfO_{2'}(X)=1$.
    Hence involutions of $\rmF$ are fused and involutions of $X$ are contained in $\rmF$ by Lemma~\ref{lem:pre}\,(1).
    Then $\rmF$ is one of $2$-groups given in Theorem~\ref{thm:p-group}.
    If $X=\rmF$, then $G\cong A{.}(B{:}C)$ with $B=\rmF$ and $C=1$ as desired.
    
    Assume that $X>\rmF$ and $\rmF$ is homocyclic.
    Then each minimal normal subgroup of $X/\rmF=X/\bfO_2(X)$ is an elementary abelian $p$-group for some odd prime $p$.
    Let $p$ be prime divisor of $|\soc(X/\rmF)|$, and let $S\lhd X$ such that $\rmF\lhd S$ and $S/\rmF\neq 1$ is the product of all minimal normal $p$-subgroups of $X/\rmF$.
    Then $S/\rmF$ is an elementary abelian $p$-group.
    This yields that $S=\rmF{:}H$, where $H\cong S/\rmF$ is an elementary abelian $p$-group.
    Note that $\bfC_{X}(\rmF)\leqslant\rmF$ and $\bfC_X(\rmF)$ is characteristic in $X$. Then $\bfC_X(\rmF)=\rmF$, and so $\bfC_H(\rmF)=1$. 
    Note that $\rmF$ is a homocyclic $2$-group and $|H|$ is a $p$-group for some odd prime $p$.
    We have that $\rmF=[\rmF,H]\times \bfC_{\rmF}(H)$, see~\cite[24.6]{aschbacher2000Finite} for instance.
    Then $\bfC_H(\rmF)=1$ implies $\bfC_{\rmF}(H)=1$.
    By Lemma~\ref{lem:splitnormal}, $X=\rmF{:}\bfN_{X}(H)$.
    Since all involutions of $X$ lie in $\rmF$, the complement $\bfN_{X}(H)$ has odd order.
    Thus $G=A.X=A.(B{:}C)=(A{:}B).C$ with $A=\bfO_{2'}(G)$, $B\cong \rmF$ and $C\cong \bfN_{X}(H)$.

    Assume that $X>\rmF$ and $\rmF$ is nonabelian.
    Theorem~\ref{thm:p-group} implies that either $\rmF\cong \rmQ_{2^n}$, or $\rmF$ is a Suzuki $2$-group.
    Then $\Aut(\rmF)$ is solvable.
    Since $\rmF$ is a characteristic subgroup of $X$, there is a natural homomorphism $\varphi:\Aut(X)\rightarrow \Aut(\rmF)$.
    We claim that the kernel of $\varphi$ is an abelian group.
    
    Indeed, let $f\in \ker\varphi$.
    Then $f$ fixes every element in $\rmF$.
    Let $\alpha: X\rightarrow X$ be a function such that $\alpha(x)=f(x)x^{-1}\in X$.
    Then $f(x)=\alpha(x) x$, and for any $z\in\rmF$, we have that 
    \[x^{-1}z x=f(x^{-1}z x)=f(x)^{-1}f(z)f(x)=x^{-1}\alpha(x)^{-1}z\alpha(x)x.\]
    This implies that $z=z^{\alpha(x)}$ for all $z\in \rmF$, and hence $\alpha(x)\in\bfC_{X}(\rmF)$.
    Since $X$ is solvable and $\rmF$ is the Fitting subgroup of $X$, we have that $\alpha(x)\in\bfC_{\rmF}(\rmF)=\bfZ(\rmF)$.
    Thus $\alpha$ is a function from $X$ to $\bfZ(\rmF)$.
    Now assume that $g\in\ker\varphi$ with $\beta: X\rightarrow X$ such that $g(x)=\beta(x)x$ for $x\in X$.
    We have that
    \[f\circ g(x)=f(\beta(x)x)=f(\beta(x))f(x)=\beta(x)\alpha(x)x, \mbox{ for $x\in X$.}\]
    Similarly, we have that $g\circ f(x)=\alpha(x)\beta(x)x$ for $x\in X$.
    Note that $\alpha(x),\beta(x)\in \bfZ(\rmF)$ for $x\in X$.
    This implies that $f\circ g=g\circ f$, and hence $\ker\varphi$ is abelian.
    The claim holds.
    
    Since both $\Aut(\rmF)$ and $\ker\varphi$ are solvable, we have that $\Aut(X)$ is solvable.
    Let $O$ be a Hall $2'$-group of $\Aut(X)$.
    Note that the number of involutions is odd and $\Aut(X)$ is transitive on the involutions of $X$.
    Then $X{:}O$ is a solvable group with a unique conjugacy class of involutions.
    By Thompson's theorem, $X{:}O$ has $2$-length $1$, so does $X$.
    Thus $X=\rmF{:} H$ for some subgroup $H<X$ of odd order.
    Therefore, we have that 
    \[G=A{.}X=A.(B{:}C)=(A{:}B){.}C,\]
    where both $A$ and $C$ have odd order and $B$ is a $2$-group given in Theorem~\ref{thm:p-group}.
\end{proof}

\subsection{Proof of Theorems~\ref{thm:main_inv_fuse} and~\ref{thm:KTS}}\

Theorem~\ref{thm:reduction} divides the proof of Theorem~\ref{thm:main_inv_fuse} into three cases, and the first two cases have already been completed: solvable group case (Theorem~\ref{thm:sol_inv_fuse}) and almost quasisimple group case (Theorem~\ref{thm:as_inv_fuse} and Lemma~\ref{lem:aqs_inv}).
In this section, we will prove that groups in the third case of Theorem~\ref{thm:reduction} are solvable, and then complete the proof of Theorem~\ref{thm:main_inv_fuse}.
We need several lemmas.
\begin{lemma}\label{lem:new_inv}
    Let $G$ be a finite group and $\bfO_2(G)\neq 1$.
    For $H<G$, denote by $\overline{H}$ the quotient image of $H$ in $G/\bfO_2(G)$. 
    If there exists $g\in G$ of odd order such that 
    \[\bfC_{\bfO_2(G)}(g)=1\mbox{ and }|\bfN_{\overline{G}}\bigl(\overline{\langle g\rangle}\bigr)|\mbox{ is even,}\]
    then $G\setminus\bfO_2(G)$ contains some involutions.
\end{lemma}
\begin{proof}
    Let $N=\bfO_2(G)$, and let $x\in G$ such that $\overline{\langle x\rangle}< \bfN_{\overline{G}}\bigl(\overline{\langle g\rangle}\bigr)$  has order $2$.
    Then $g^x=vg^\ell$ for some $v\in N$ and integer $\ell$.
    Note that 
    \[\langle N,g^x\rangle=N{:}\langle g^x\rangle=N{:}\langle vg^\ell\rangle=N{:}\langle g\rangle.\]
    By Schur-Zassenhaus' Theorem, there exists $w\in N$ such that $\langle g^x\rangle ^w=\langle g\rangle$.
    Then $g^{xw}\in\langle g\rangle$, and hence $xw\in \bfN_G(\langle g\rangle)$.
    Since $\bfN_G(\langle g\rangle)\cap N=\bfN_N(\langle g\rangle)=\bfC_N(g)=1$, we have that $\langle xw\rangle\cap N=1$.
    Note that $|\overline{\langle xw\rangle}|=|\overline{\langle x\rangle}|$ is even.
    The group $\langle xw\rangle$ is of even order and contains an involution $y$.
    Thus $y\in G\setminus N$ is an involution, which completes the proof.
\end{proof}

\begin{lemma}\label{lem:trans_fix_free}
    Let $G\leqslant \GL_n(2)$ be a nonsolvable transitive linear group acting on $V=\bbF_2^n$.
    Then either 
    \begin{enumerate}
        \item $G$ contains an element $g$ of odd order such that the only fixed point of $g$ is the zero vector and $\bfN_G(g)$ has even order; or
        \item $\SL_3(2^f)\leqslant G\leqslant \GammaL_3(2^f)$ and $V=\bbF_2^{3f}$.
    \end{enumerate}
\end{lemma}
\begin{proof}
    By the classification of finite transitive linear groups, a nonsolvable transitive linear group $G$ on $V=\bbF_2^n$ satisfies one of the following cases.
    \begin{enumerate}
        \item[1.] $\SL_m(2^f)\trianglelefteqslant G\leqslant\GammaL_m(2^f)$ with $n=mf$;
        \item[2.] $\Sp_{2a}(2^f)\trianglelefteqslant G\leqslant \CGammaSp_{2a}(2^f)$ with $n=2af$;
        \item[3.] $\rmG_2(2^f)'\trianglelefteqslant G\leqslant \rmG_2(2^f){:}\rmC_f$ with $n=6f$;
        \item[4.] $G\cong \rmA_6$ or $\rmA_7$ for $V=\bbF_2^4$;
    \end{enumerate}
    Note that if $H<G$ are both transitive linear groups and $H$ satisfies part~(1), then $G$ satisfies part~(1).

    We claim that if $G$ lies in Cases~1 and 2 and contains a $\calC_3$-subgroup $H\cong \SL_2(2^f)$ with $n=2f$, then $G$ satisfies part~(1).
    Let $H$ be a $\calC_3$-subgroup isomorphic to $\SL_2(2^f)$. 
    Note that $H$ contains a Singer cycle $L\cong\rmC_{2^{k}+1}$ and the only fixed point of $L$ is the zero vector.
    Let $g$ be a generator of $L$.
    Then  $\bfN_G(L)\geqslant \bfN_H(L)\cong \rmD_{2(2^k+1)}$ has even order.
    Thus $G$ satisfies part~(1).

    Now we show that if $G$ lies in Cases~1 and 2, then one of parts~(1) and (2) holds.
    Suppose that $G$ does not satisfy part~(1).
    By Step 1, we have that $\SL_{2d+1}(2^f)\trianglelefteqslant G\leqslant \GammaL_{2d+1}(2^f)$ with $n=(2d+1)f$ and $d\geqslant 1$.
    If $d=1$, then part~(2) holds.
    Suppose that $d>1$.
    Then $V$ has a decomposition $V=V_0\oplus V_1\oplus\cdots \oplus V_{d-1}$, where $\dim V_0=3$ and $\dim V_i=2$ for $i=1,...,d-1$.
    Hence $G$ contains a subgroup 
    \[H=\SL(V_0)\times \SL(V_1)\times\cdots\times \SL(V_{d-1})\cong \SL_3(2^f)\times \SL_2(2^f)^{d-1}.\]
    Let $L_i$ be a Singer cycle in $\SL(V_i)$, and let $g_i$ be a generator of $L_i$ for $i=0,...,d-1$.
    Set $g=g_0g_1\cdots g_{d-1}$.
    Then the zero vector is the unique fixed point of $g$.
    Note that 
    \[\bfN_G(\langle g\rangle)\geqslant \bfN_H(\langle g\rangle)\gtrsim \bfN_{\SL(V_0)}(\langle g_0\rangle)\times \rmD_{2(2^f+1)}\cong (\rmC_{2^{2f}+2^f+1}{:}\rmC_3)\times \rmD_{2(2^f+1)}.\]
    Hence $G$ satisfies part~(1).
    
    Suppose that $G=\rmG_2(2)'\cong \PSU_3(3)$ for $V=\bbF_2^6$ or $G\cong \rmA_6$ or $\rmA_7$ for $V=\bbF_2^4$.
    With the aid of \Magma, we have that $G$ satisfies part~(1).
    Now we only need to deal with the case where $\rmG_2(2^{f})\trianglelefteqslant G\leqslant \rmG_2(2^f){:}\rmC_f$ with $n=6f$.
    Note that $G$ contains a subgroup $H$ isomorphic to $\rmG_2(2)$ such that $V=V_1\oplus \cdots \oplus V_f$, where $V_i\cong\bbF_2^6$ and $V_i\cong V_j$ as $\bbF_2H$-modules for $1\leqslant i,j\leqslant f$.
    Note $H'\cong \rmG_2(2)'\cong \PSU_3(3)$ contains $g$ such that $|g|$ is odd with no fixed point on $V_i\setminus\{0\}$ for each $i$ and $\bfN_{H'}(\langle g\rangle)$ has even order.
    Thus $\bfN_G(\langle g\rangle)$ has even order and $g$ has no fixed point on $V\setminus\{0\}$.
    Then part~(1) holds.
    Since we have exhausted all possibilities of $G$, the proof is complete.
\end{proof}

Now we are ready to show that groups in third part of Theorem~\ref{thm:reduction} are solvable.
\begin{lemma}\label{lem:case_3_sol}
    Assume that $G$ has a unique conjugacy class of involutions such that $\bfO_2(G)\cong \rmC_{2^m}^n$ and $G/\bfO_2(G)\lesssim \GL_n(2)$ is faithful and transitive on non-identities of both 
    \[\Omega_1(\bfO_2(G))=\{x\in\bfO_2(G):x^2=1\}\mbox{ and }\bfO_2(G)/\Phi(\bfO_2(G)).\]
    Then $G$ is solvable.
\end{lemma}
\begin{proof}
    Suppose that $G$ is nonsolvable.
    Let $N=\bfO_2(G)\cong \rmC_{2^m}^n$ with $m,n\geqslant 1$.
    Then $G/N\lesssim\GL_n(2)$.
    Suppose that $\overline{g}\in G/N$ such that $|\overline{g}|$ is odd, $\overline{g}$ has no fixed point on $\Omega_1(N)\setminus\{1\}$ and $\bfN_{G/N}(\langle \overline{g}\rangle)$ has even order.
    Then there exists $g\in G$ such that $gN=\overline{g}$ and $g$ has odd order.
    This implies that $\bfC_N(g)\cap \Omega_1(N)=1$, and then $\bfC_N(g)=1$.
    Thus $g$ has no fixed points on $N\setminus\{1\}$, so does $\overline{g}$.
    By Lemma~\ref{lem:new_inv}, $G\setminus N$ contains some involutions, and thus $G$ has more than one conjugacy class of involutions, a contradiction.
    Thus $(G/N)^{(\infty)}\cong \SL_3(2^f)$ with $n=3f$ by Lemma~\ref{lem:trans_fix_free}.

    Note that $G/N$ naturally contains a subgroup $\overline{H}$ isomorphic to $\SL_3(2)$ such that $N/\Phi(N)$ can be viewed as a direct sum of $f$ copies of the natural $\overline{H}$-module $W\cong\bbF_2^3$.
    The stabilizer in $\overline{H}$ of an nonzero vector of $W$ is isomorphic to $C_2^2{:}\SL_2(2)\cong \rmS_4$.
    Then this stabilizer contains a unique subgroup $\overline{L}$ isomorphic to $\rmA_4$.
    Let $L$ be the full preimage of $\overline{L}$ in $G$.
    Then 
    \[L=N.L\cong \rmC_{2^m}^{3f}.\rmA_4.\] 
    With the aid of \Magma, we have that $\rmH^2(\overline{L},W)=0$.
    This implies that $L/\Phi(N)\cong \rmC_{2}^{3f}{:}\rmA_4$ is a split extension.
    Then there exists $L_1<L$ such that $\Phi(N)\lhd L_1$ and 
    \[L/\Phi(N)= (N/\Phi(N)){:}(L_1/\Phi(N))\cong \rmC_2^{3f}{:}\rmA_4.\]
    Note that the action of $L_1/\Phi(N)\cong \rmA_4$ on $\Phi(N)/\Phi(\Phi(N))$ is isomorphic to the action of $L/N\cong\rmA_4$ on $N/\Phi(N)$.
    Thus the previous process can be repeated, and finally we can obtain a subgroup $L_k$ of $L<G$ isomorphic to $\rmA_4$ and $L_k\cap N=1$.
    Since $L_k$ has involutions, $G$ has some involutions that are not in $N=\bfO_2(G)$.
    This implies that $G$ has more than one conjugacy class of involutions, a contradiction.
    Therefore, $G$ is solvable.
\end{proof}

We summarize the proof of Theorem~\ref{thm:main_inv_fuse} as below.

\begin{proof}[Proof of Theorem~\ref{thm:main_inv_fuse}]
    If $G$ is solvable, then Theorem~\ref{thm:sol_inv_fuse} implies that $G$ has $2$-length $1$ as in part~(1).
    Assume that $G$ is nonsolvable.
    Let $\overline{G}=G/\bfO_{2'}(G)$.
    Then $\overline{G}$ is nonsolvable and has a unique fusion class of involutions by Lemma~\ref{lem:pre}\,(2).
    By Theorem~\ref{thm:reduction} and Lemma~\ref{lem:case_3_sol}, we deduce that $\overline{G}$ is an almost quasisimple group.
    If $\bfO_2(\overline{G})\neq 1$, then $\overline{G}$ has a unique involution and $\overline{G}^{(\infty)}\cong \SL_2(q)$ or $2.\rmA_7$ by Lemma~\ref{lem:aqs_inv}, as in (i) of part~(2).
    Assume that $\bfO_2(\overline{G})=1$.
    Then $\soc(\overline{G})$ is simple as $\bfO_{2'}(\overline{G})=1$, and hence $\overline{G}$ is almost simple.
    We have that $\overline{G}$ has a unique conjugacy class of involutions, and $\overline{G}^{(\infty)}$ is one of the groups in Table~\ref{tab:pyra}.
    Let $x$ and $y$ be two involutions in $G$, and let $\overline{x}$ and $\overline{y}$ be quotient images of $x$ and $y$ in $\overline{G}$, respectively.
    Then there exists $z\in G$ such that $\overline{x^z}=\overline{y}$, and hence 
    \[x^z\in \langle \bfO_{2'}(G),y\rangle=\bfO_{2'}(G){:}\langle y\rangle.\]
    By Schur-Zassenhaus' Theorem, involutions in $\langle \bfO_{2'}(G),y\rangle$ are conjugate.
    Thus there exists $z_0\in G$ such that $(x^z)^{z_0}=y$, and therefore involutions of $G$ are conjugate, as in part~(2).
    The proof is complete.
\end{proof}

Recall that an $m$-pyramidal group is a finite group which has precisely $m$ involutions and these $m$ involutions are pairwise conjugate.
Before proving Theorem~\ref{thm:KTS}, we give the following lemma.
\begin{lemma}\label{lem:pyra_quo}
    Assume that $G$ is an $m$-pyramidal group with $G'$ of even order.
    Let $O\lhd G$ be a normal subgroup of odd order, and let $x$ be an involution.
    Then $G/O$ is a $k$-pyramidal group, and $m/k=|O|/|\bfC_O(x)|$ is an odd square number.
\end{lemma}
\begin{proof}
    As $O$ is of odd order, we can take $\{g_1O,...,g_kO\}$ as the set of involutions in $G/O$ such that each $g_i$ is an involution in $G$.
    Hence each involution in $G$ has the form $g_i w$ for some $w\in O$ and $i=1,...,k$.
    Let $\Omega_i=\{w\in O: w^{g_i}=w^{-1}\}$.
    Then $|\Omega_i|=|O|/|\bfC_O(g_i)|$ and, for $w\in O$, $g_iw$ is an involution if and only if $w\in\Omega_i$.
    As $g_i$ is conjugate to $g_j$ for each $i,j$, we have that $|\Omega_i|=|\Omega_j|=|O|/|\bfC_O(x)|$.
    Thus $m=k|\Omega_i|=k|O|/|\bfC_O(x)|$.

    We finish the proof by claiming that $|\Omega_i|$ is a square number.
    Assume that the claim is false and choose a counterexample $(G,O,x)$ with $|O|$ minimal.
    Then $O\neq 1$.
    Suppose that $O$ is not a minimal normal subgroup.
    Let $N$ be a minimal normal subgroup of $G$ contained in $O$.
    Then one of $(G,N,x)$ and $(G/N,O/N,xN)$ is a counterexample of the lemma, and both $|N|$ and $|O/N|$ are smaller than $|O|$, a contradiction.
    Thus $O$ is a minimal normal subgroup, and hence $O\cong \rmC_p^n$ for some odd prime $p$.
    Then there is a natural group homomorphism $\varphi:G\rightarrow \GL_n(p)$ induced by conjugating on $O$.
    Note that $|\Omega_i|$ equals to the size of eigenspace of $\varphi(x)$ with respect to eigenvalue $-1$.
    As $G'$ contains an involution, we have that $\varphi(x)$ must lie in $\SL_n(p)$, and hence has determinant $1$.
    Thus the eigenspace of $\varphi(x)$ with respect to eigenvalue $-1$ has even dimension.
    This implies that $|\Omega_i|$ is a square number, a contradiction.
    Therefore, the proof is complete.
\end{proof}

We finish the paper by proving Theorem~\ref{thm:KTS}.
\begin{proof}[Proof of Theorem~\ref{thm:KTS}]
    By \cite[Proposition~1.1]{gao2025}, $G$ is an $m$-pyramidal group.
    Since $G$ is nonsolvable, $G/O$ is an almost quasisimple $k$-pyramidal group by Theorem~\ref{thm:main_inv_fuse}, and hence part~(1) holds.
    Note that $G'$ is of even order.
    By Lemma~\ref{lem:pyra_quo}, $m/k=|O|/|\bfC_O(x)|$ is an odd square number for each involution $x\in G$, as in part~(2).

    Now we prove part~(3).
    If $3\mid m/k$, then case~(i) holds and there is nothing to prove.
    Hence we assume that $m/k=|O|/|\bfC_O(x)|$ is not divisible by $3$.
    Note that $\KTS(v)$ exists if and only if $v\equiv 3\pmod 6$, see~\cite{Ray1971}.
    Then an $m$-pyramidal $\KTS(v)$ realized under $G$ exists only if $|G|+m\equiv 3\pmod 6$.
    Since $G$ is nonsolvable, it has even order.
    Note that the number of involutions of an even order group is odd.
    This implies that $|G|+m$ is always odd, and hence an $m$-pyramidal $\KTS(v)$ realized under $G$ exists only if $|G|+m$ is divisible by $3$.
    Suppose that $|G|$ is divisible by $3$.
    Since $m/k$ is an odd square number, $m/k\equiv 1\pmod 3$.
    Then 
    \[|G|+m\equiv m\equiv k \pmod 3.\]
    Hence if $3\nmid m/k$ and $|G|$ is divisible by $3$, such $\KTS(v)$ exists only if $3\mid k$, and thus case (ii) holds.
    Suppose that $|G|$ is not divisible by $3$.
    Then the order of the quasisimple $(G/O)^{(\infty)}$ is not divisible by $3$.
    This implies that $(G/O)^{(\infty)}\cong {}^2\rmB_2(q)$ with $q=2^{2n+1}$ for $n\geqslant 1$ and $O$ is a $3'$-group.
    Since $O$ has odd order, it is a $\{2,3\}'$-group, and then case~(iii) holds.
    Thus part~(3) holds.

    Now assume that $m$ is a prime power and $G$ is a nonsolvable $m$-pyramidal group.
    Since the number of involutions of an even order group is odd, $m$ is a power of an odd prime.
    It follows that either $k=1$ or $k$ is also a power of an odd prime.
    It is routine to check from the displayed factorizations of Table~\ref{tab:pyra} that, except for the first two rows $k=1$, the number $k$ is composite.
    Thus $G/O$ has a unique involution and $m$ is a square of some odd number by part~(2), and hence $m=p^{2f}$ for some $f\geqslant 1$.
    Suppose that $m=9$.
    We may suppose that $G$ is a nonsolvable $9$-pyramidal group of the minimal order.
    Let $N$ be a minimal normal subgroup of $G$.
    Then $N$ is elementary abelian as there is no simple group with $9$ involutions and all its involutions are fused by Theorem~\ref{thm:main_inv_fuse}.
    Thus $G/N$ has a unique involution by the minimality of $|G|$.
    Let $x$ be an involution of $G$.
    Then, for any $g\in G$, $x^g=zx$ for some $z\in N$.
    Thus 
    \[\bfC_N(x)^g=\bfC_N(zx)=\bfC_N(x).\]
    This implies that $\bfC_N(x)$ is a normal subgroup of $G$.
    It follows that $G/\bfC_N(x)$ is a nonsolvable $9$-pyramidal group, contradicting the minimality of $G$.
    Hence we have that $\bfC_N(x)=1$, and then $|N|=m|\bfC_N(x)|=9$.
    Thus $N\cong\rmC_3^2$.
    Since $x\notin\bfC_{X}(N)$, $X/\bfC_X(N)$ contains a nonabelian composition factor isomorphic to $\rmA_7$ or $\PSL_2(q)$.
    This is impossible as $X/\bfC_X(N)\lesssim\Aut(N)\cong\GL_2(3)$ is solvable.
    Thus there are no nonsolvable $9$-pyramidal groups.
    Note that involutions of $\ASL_2(p^f)\cong p^{2f}{:}\SL_2(p^f)$ are conjugate and $\ASL_2(p^f)$ is nonsolvable when $p^{2f}>9$.
    It follows that there exist nonsolvable $m$-pyramidal groups for $m=p^{2f}>9$.
    Recall that $m$ is always odd as the number of involutions is always odd.
    Hence, for a prime power $m$, there exists a nonsolvable $m$-pyramidal group if and only if $m=p^{2f}>9$ with $p$ odd prime. 
    Therefore, the ``in particular'' part holds.
\end{proof}

\vskip 3mm
\noindent{\bf AI declaration}:
The main results of this paper were obtained in March 2026. 
No AI tool was used in the mathematical work, including the derivation and mathematical formulation of the results.
During manuscript preparation, the authors used DeepSeek V4.1 Flash and a later version (1017), together with DSH (DeepSeek Harness), only for literature search and language polishing. 
These tools did not generate mathematical content or change the mathematical substance of the paper. 
Language polishing was limited to grammar, wording, and readability and did not affect the mathematical statements or proofs. 
All references were checked by the authors, who take full responsibility for the article.

\vskip 3mm

\noindent{\bf Acknowledgement}:
This work was partially supported by NSFC: 12601042;
Foundation for Cultivated Young Talents of Fujian Province, China: 202501AT070016;  
Fundamental Research Funds for the Central Universities: 20720260155.

\end{document}